\documentclass[11pt]{article}

\usepackage{amssymb,latexsym,amsmath,amsthm, mathrsfs}
\usepackage{graphicx}
\usepackage[all]{xy}
\newtheorem{theorem}{Theorem}[section]

\newtheorem{proposition}[theorem]{Proposition}
\newtheorem{corollary}[theorem]{Corollary}

\newenvironment{definition}[1][Definition]{\begin{trivlist}
\item[\hskip \labelsep {\bfseries #1}]}{\end{trivlist}}
\newenvironment{example}[1][Example]{\begin{trivlist}
\item[\hskip \labelsep {\bfseries #1}]}{\end{trivlist}}

\newenvironment{conjecture}[1][Conjecture]
{\begin{trivlist} \item[\hskip \labelsep {\bfseries #1}]}
{\end{trivlist}}
\newenvironment{acknowledgement}[1][Acknowledgements]
{\begin{trivlist} \item[\hskip \labelsep {\bfseries #1}]}
{\end{trivlist}}

\DeclareMathSymbol{\N}{\mathbin}{AMSb}{"4E}
\DeclareMathSymbol{\Z}{\mathbin}{AMSb}{"5A}
\DeclareMathSymbol{\R}{\mathbin}{AMSb}{"52}
\DeclareMathSymbol{\Q}{\mathbin}{AMSb}{"51}
\DeclareMathSymbol{\I}{\mathbin}{AMSb}{"49}
\DeclareMathSymbol{\C}{\mathbin}{AMSb}{"43}

\def\P{{\mathbb P}}

\def\A{{\mathbb A}}
\def\F{{\mathbb F}}

\title{Grothendieck standard conjectures, morphic cohomology and Hodge index theorem}
\author{Jyh-Haur Teh}

\begin{document}
\maketitle

\footnotetext[1]{Tel: +886 3 5712121-33071; Fax: +886 3 5723888\\
E-mail address: jyhhaur@math.nthu.edu.tw}

\begin{abstract}
Using morphic cohomology, we produce a sequence of conjectures,
called morphic conjectures, which terminates at the Grothendieck
standard conjecture A. A refinement of Hodge structures is given,
and with the assumption of morphic conjectures, we prove a Hodge
index theorem. We answer a question of Friedlander and Lawson by
assuming the Grothendieck standard conjecture B and prove that the
topological filtration from morphic cohomology is equal to the
Grothendieck arithmetic filtration for some cases.
\end{abstract}
\section{Introduction}
The homotopy groups of the cycle spaces of a complex projective
variety $X$ form a set of invariants, called the Lawson homology
groups of $X$(see \cite{F1}, \cite{L1}). To establish a
cohomology-like theory, Friedlander and Lawson produced the notion
of algebraic cocycle in \cite{FL1} and defined the morphic
cohomology groups of a projective variety to be the homotopy groups
of some algebraic cocycle spaces. Furthermore, for a smooth
projective variety $X$, by using their moving lemma (see
\cite{FL2}), they proved a duality theorem between the Lawson
homology and morphic cohomology of $X$. Walker has defined an
inductive limit of mixed Hodge structures on the Lawson homology
\cite{W} of $X$, and we extend this to morphic cohomology by using
the above duality isomorphisms. Then the images of the morphic
cohomology groups of $X$ in its singular cohomology groups under the
natural transformations have sub-Hodge structures.

The Grothendieck standard conjectures have various parts (see
\cite{Groth}, \cite{Kle1}, \cite{Kle2}). For a smooth projective
variety $X$ of dimension $m$, let $C^j(X)$ be the subspace of
$H^{2j}(X; \Q)$ which is generated by algebraic cycles. By the Hard
Lefschetz theorem, cup product with the Lefschetz class $L$ gives
isomorphism
$$\xymatrix{ H^{2j}(X; \Q) \ar[r]^{L^{m-2j}}_{\cong} & H^{2m-2j}(X;
\Q )\\
C^j(X) \ar[u] \ar[r] & C^{m-j}(x) \ar[u]}$$ for $j\leq \lfloor
\frac{m}{2} \rfloor$. The Grothendieck standard conjecture A (GSCA
for short) claims that the restriction of $L^{m-2j}$ also gives an
isomorphism between $C^j(X)$ and $C^{m-j}(X)$, or equivalently, the
adjoint operator $\Lambda$ maps $C^{m-j}(X)$ into $C^j(X)$. The
Grothendieck standard conjecture B (GSCB for short) says that the
adjoint operator $\Lambda$ is algebraic, i.e., there is a cycle
$\beta$ on $X\times X$ such that $\Lambda: H^*(X; \Q)
\longrightarrow H^*(X; \Q)$ is got by lifting a class from $X$ to
$X\times X$ by the first projection, cupping with $\beta$ and taking
the image in $H^*(X; \Q)$ by the Gysin homomorphism associated to
the second projection. For abelian varieties, the GSCB was proved by
Lieberman in \cite{Lieb}, and we know the GSCB for a smooth variety
which is a complete intersection in some projective space and for
Grassmannians (see \cite{Groth}).

In this paper, a sequence of conjectures, called morphic
conjectures, is introduced and the GSCA is the last conjecture in
this sequence. We show that if the GSCB holds on $X$, it implies all
the morphic conjectures of $X$. Various equivalent forms of the
morphic conjectures are provided. It is well known that the GSCA is
equivalent to the statement that numerical equivalence is equal to
homological equivalence. We prove an analogous statement for our
morphic conjectures in Proposition \ref{equivalent forms}. It was
proved by Jannsen (see \cite{Jannsen}) that the GSCA is equivalent
to the semisimplicity of the ring of algebraic correspondences, we
do not know if analogous result is true for morphic conjectures. The
refinement of the Hodge structures by the images of the morphic
cohomology groups of $X$ in its singular cohomology groups, with the
assumption of the corresponding morphic conjecture, is compatible
with a refinement of the Lefschetz decomposition, and we get a
result analogous to the classical Hodge index theorem.

Let us give a brief review of morphic cohomology (see \cite{FL3}).
Throughout this paper, $X, Y$ are smooth complex projective
varieties and the dimension of $X$ is $m$. An effective $Y$-valued
$r$-cocycle $c$ is an effective $(r+m)$-cycle in $X\times Y$ such
that each fibre of $c$ over $X$ is an $r$-cycle on $Y$. We denote
the group of effective $Y$-valued $r$-cocycles by
$\mathscr{C}_r(Y)(X)$. The Chow monoid $\mathscr{C}_{r+m}(X\times
Y):=\coprod_{d\geq 0}\mathscr{C}_{r+m, d}(X\times Y)$ with topology
from the analytic topology of each Chow variety $\mathscr{C}_{r+m,
d}(X\times Y)$ is a topological monoid and we give
$\mathscr{C}_r(Y)(X) \hookrightarrow \mathscr{C}_{r+m}(X\times Y)$
the subspace topology. Let $Z_r(Y)(X):=\mathscr{C}_r(Y)(X)\times
\mathscr{C}_r(Y)(X)/\sim$ be the naive group-completion of
$\mathscr{C}_r(Y)(X)$ with the quotient topology where $(a, b)\sim
(c, d)$ if and only if $a+d=b+c$, then $Z_r(Y)(X)$ is a topological
abelian group. Define the group of algebraic $t$-cocycles on $X$ to
be
$$Z^t(X):=\frac{Z_0(\P^t)(X)}{Z_0(\P^{t-1})(X)}$$
and define the $(t, k)$-morphic cohomology group to be
$$L^tH^k(X):=\pi_{2t-k}Z^t(X), \ \ 0\leq k \leq 2t$$ the homotopy group of the cocycle
space. It is shown in \cite[Theorem A.2]{T1} that $Z_0(\P^t)(X)$ is
a CW-complex, and a similar argument of \cite[Theorem A.5]{T1} shows
that $Z^t(X)$ is a CW-complex. We list some fundamental properties
of morphic cohomology which will be used in this paper.

\begin{enumerate}
\item There is a natural transformation $\Phi^{t, k}: L^tH^k(X)
\longrightarrow H^k(X)$ from morphic cohomology to singular
cohomology for any $t, k$.

\item There is a cup product pairing $L^tH^k(X)\otimes L^rH^s(X)
\longrightarrow L^{t+r}H^{k+s}(X)$ which is transformed to the cup
product in singular cohomology by the natural transformations.

\item There is a commutative diagram:
$$\xymatrix{ L^tH^k(X) \ar[r]^-{\mathscr{D}} \ar[d]_{\Phi} & L_{m-t}H_{2m-k}(X) \ar[d]^{\Psi}\\
H^k(X) \ar[r]^-{\mathscr{PD}} & H_{2m-k}(X)}
$$ for $0\leq t \leq m$ where $\mathscr{D}$ is the Friedlander-Lawson duality
isomorphism, $\mathscr{PD}$ is the Poincar\'e duality isomorphism,
and $\Psi$ is the natural transformation from Lawson homology to
singular homology.
\end{enumerate}

In \cite{FL1}, Question 9.7, Lawson and Friedlander asked if the map
$\Phi^{t, k}:L^tH^k(X) \longrightarrow H^k(X)$ is surjective for $X$
smooth and $t\geq k$. We answer this question in rational
coefficients by assuming GSCB.  As a consequence, we prove that the
topological filtration from morphic cohomology is equal to the
Grothendieck arithmetic filtration for some cases. This may open a
new way to check the validity of the Grothendieck standard
conjectures and the generalized Hodge conjecture.

\section{Inductive limit of mixed Hodge structures}
We use HS and MHS to abbreviate  Hodge structure and mixed Hodge
structure respectively. We follow Walker's definition of the
inductive limit of mixed Hodge structures (IMHS) in \cite{W}.

\begin{definition}
An IMHS is an inductive system of MHS's $\{H_{\alpha}, \alpha\in
I\}$ where the index set $I$ is countable such that there exist
integers $M<N$ so that $W_M((H_{\alpha})_{\Q})=0,
W_N((H_{\alpha})_{\Q})=(H_{\alpha})_{\Q}, F^N((H_{\alpha})_{\Q})=0,
\mbox{ and } F^M((H_{\alpha})_{\Q})=(H_{\alpha})_{\Q}$ for all
$\alpha\in I$. Equivalently, an IMHS is a triple $(H, W_{\bullet},
F^{\bullet})$, where $H$ is a countable abelian group,
$W_{\bullet}(H_{\Q})$ and $F^{\bullet}(H_{\C})$ are finite
filtrations satisfying
$$Gr^W_n(H_{\C})=\underset{p+q=n}\oplus H^{p, q}$$ where
$$H^{p, q}=F^pGr^W_{p+q}(H_{\C})\cap
\overline{F}^qGr^W_{p+q}(H_{\C}),$$ and such that every finitely
generated subgroup of $H$ is contained in a finitely generated
subgroup $H'$ so that $(H', W_{\bullet}|_{H'_{\Q}},
F^{\bullet}|_{H'_{\C}})$ is a MHS. A morphism of IMHS is morphism of
filtered systems of MHS's.
\end{definition}

It is shown in Proposition 1.4 of \cite{Hu} that the category of
IMHS is abelian. By \cite[Theorem 4.1]{W}, the Lawson homology
groups of a quasi-projective variety have an IMHS.

\begin{definition}
We endow $L^tH^k(X)$ with an IMHS by making the Friedlander-Lawson
duality map $\mathscr{D}:L^tH^k(X) \longrightarrow
L_{m-t}H_{2m-k}(X)$ an isomorphism of IMHS.
\end{definition}

\begin{proposition}
The map $\Phi^{t, k}: L^tH^k(X) \longrightarrow H^k(X)$ is a
morphism of IMHS and the IMHS on Im$\Phi^{t, k}$ is a sub-HS of
$H^k(X)$. \label{natural transformation morphism}
\end{proposition}

\begin{proof}
Consider the following commutative diagram:
$$\xymatrix{ L^tH^k(X) \ar[r]^-{\mathscr{D}} \ar[d]_{\Phi^{t, k}} &
L_{m-t}H_{2m-k}(X) \ar[d]^{\Psi_{m-t, 2m-k}}\\
H^k(X) \ar[r]^-{\mathscr{PD}} & H_{2m-k}(X)}
$$ By \cite[Theorem 4.1]{W}, $\Psi_{m-k, 2m-k}$ is
a morphism of IMHS. Since $\mathscr{PD}^{-1}$ is a morphism of HS,
$\Phi^{t, k}=\mathscr{PD}^{-1}\circ \Psi_{m-t, 2m-k} \circ
\mathscr{D}$ is a morphism of IMHS. The IMHS on $H^k(X)$ is a HS,
therefore the image of $\Phi^{t, k}$ is a sub-HS of $H^k(X)$.
\end{proof}

\begin{definition}
Define $\widetilde{L}^tH^k(X)=Im\Phi^{t, k}$ and decompose
$$\widetilde{L}^tH^k(X; \C)=\bigoplus_{p+q=k}H^{p, q}_t(X).$$
We define the morphic Hodge numbers of $X$ to be
$$h^{p, q}_t(X)=dimH^{p, q}_t(X).$$
\end{definition}

The following result (see \cite{FL1}, Theorem 4.4) says that the
image of $\Phi$ is contained in a specific range in the Hodge
decomposition.

\begin{theorem}
For a smooth projective variety $X$ there is an inclusion
$$\widetilde{L}^tH^k(X; \C)\subset
\bigoplus_{\substack{p+q=k \\ |p-q|\leq 2t-k}}H^{p, q}(X)$$
\end{theorem}

The space of the most interest to us is $\widetilde{L}^tH^{2t}(X;
\Q)=H^{t, t}_t(X; \Q)=H^{t, t}_t(X)\cap H^{t, t}(X; \Q)$ which is
the space generated by algebraic cycles with rational coefficients
where $H^{t, t}(X; \Q)=H^{2t}(X; \Q)\cap H^{t, t}(X)$. We recall
that the Hodge conjecture says that $\widetilde{L}^tH^{2t}(X;
\Q)=H^{t, t}(X; \Q)$.

\section{Signatures}
Before we proceed to the definition of morphic signatures, we need
the following result in which we correct some part of the proof in
\cite[Theorem 5.8]{FL3}. We use $\F$ to indicate any one of the
fields $\Q, \R, \C$ and define $L^tH^k(X; \F)=L^tH^k(X)\otimes \F$.

Recall that there is a $S$-map (see \cite[Theorem 5.2]{FL1}) which
makes the following diagram commutes:
$$\xymatrix{L^tH^k(X; \F) \ar[rd]_{\Phi^{t, k}} \ar[rr]^S &  & \ar[ld]^{\Phi^{t+1, k}}
L^{t+1}H^k(X; \F) \\
& H^k(X; \F) & }$$ And the natural transformation $\Phi: L^*H^k(X)
\rightarrow H^k(X)$ is induced by the map $i': Z^t(X) \rightarrow
Map(X, Z_0(\C^t))$ where $Map(X, Z_0(\C^t))$ is the space of
continuous maps from $X$ to $Z_0(\C^t)$ with the compact-open
topology and $i'$ is the map induced by the inclusion map
$i:\mathscr{C}_0(\P^t)(X) \hookrightarrow Map(X,
\mathscr{C}_0(\P^t))$. Recall that $H^k(X)=\pi_{2t-k}Map(X,
Z_0(\C^t))$.

\begin{proposition}
For any $t\geq m$, the three maps in the diagram of the $S$-map are
isomorphisms for $0\leq k\leq 2m$. \label{greater than the
dimension}
\end{proposition}

\begin{proof}
Obviously it is enough to prove the statement for $\Q$-coefficients.
We have the following commutative diagrams:
$$\xymatrix{ Z_0(\P^{t-1})(X) \ar[r] \ar[d]_{\mathscr{D}} & Z_0(\P^t)(X) \ar[r] \ar[d]_{\mathscr{D}} & Z^t(X)
\ar[d]_{\mathscr{D}}\\
Z_m(X\times \P^{t-1}) \ar[r] & Z_m(X\times \P^t) \ar[r] &
Z_m(X\times \C^t)}$$ By the Friedlander-Lawson duality theorem (see
\cite{FL3}, Theorem 3.3), the first two $\mathscr{D}$ are homotopy
equivalences, and by \cite[Proposition 3.2]{T1}, the upper and lower
rows induce long exact sequences of homotopy groups
$$\xymatrix{ \cdots \ar[r] & \pi_k Z_0(\P^{t-1})(X) \ar[r] \ar[d] &
\pi_k Z_0(\P^t)(X) \ar[r] \ar[d] & \pi_k Z^t(X) \ar[r] \ar[d] &
\pi_{k-1}Z_0(\P^{t-1})(X) \ar[r] \ar[d] & \cdots \\
\cdots \ar[r] & \pi_k Z_m(X\times \P^{t-1}) \ar[r] & \pi_k
Z_m(X\times \P^t) \ar[r] & \pi_k Z_m(X\times \C^t) \ar[r] &
\pi_{k-1}Z_m(X\times \P^{t-1}) \ar[r] & \cdots }$$ By the five-lemma
we have $\mathscr{D}_*:\pi_kZ^t(X) \overset{\cong}{\longrightarrow}
\pi_kZ_m(X\times \C^k)$ for $k\geq 0$, then by the Whitehead
theorem, we know that $Z^t(X)$ is homotopy equivalent to
$Z_m(X\times \C^t)$.

If $t\geq m$, $Z_m(X\times \C^t)\cong Z_0(X\times \C^{t-m})$ by the
homotopy property of trivial bundles (see \cite[Proposition
2.3]{FG}). Applying \cite[Proposition 3.2]{T1} to the following
sequence:
$$\xymatrix{Z_0(X\times \P^{t-1}) \ar[r]^i & Z_0(X\times \P^t) \ar[r] & Z_0(X\times \A^{t-m})}$$
we get a long exact sequence of homotopy groups,
$$\cdots \longrightarrow \pi_kZ_0(X\times \P^{t-m-1})\overset{i_*}{\longrightarrow} \pi_kZ_0(X\times \P^{t-m})
\longrightarrow \pi_kZ_0(X\times \A^{m-t})\longrightarrow
\pi_{k-1}Z_0(X\times \P^{t-m-1}) \longrightarrow \cdots$$

Recall that the Dold-Thom theorem (\cite[6.10]{DT}) says that for a
CW-complex $A$,  $\pi_kZ_0(A)\cong H^{BM}_k(A; \Z)$ for all $k$.
Applying the Dold-Thom theorem and tensor by $\Q$, we get a long
exact sequence:
$$\cdots \longrightarrow H_k(X\times \P^{t-m-1};
\Q)\overset{i_*}{\longrightarrow} H_k(X\times \P^{t-m};
\Q)\longrightarrow H^{BM}_k(X\times \A^{m-t}; \Q)$$
$$\longrightarrow H_{k-1}(X\times \P^{t-m-1}; \Q) \longrightarrow
\cdots$$ where $i_*$ is induced from the inclusion map $i: X\times
\P^{t-m-1} \subset X \times \P^{t-m}$. Since the inclusion map
$j:\P^{t-m-1} \subset \P^{t-m}$ induces an isomorphism $j_*:
H_k(\P^{t-m-1}) \longrightarrow H_k(\P^{t-m})$ in homology groups
for $0\leq k\leq 2(t-m-1)$, by the K\"{u}nneth formula in homology
for $H_k(X\times \P^{t-m-1}; \Q)$ and $H_k(X\times \P^{t-m}; \Q)$,
it is not difficult to see that
$$H^{BM}_k(X\times \C^{t-m}; \Q)=
\left\{
\begin{array}{ll}
    0, & \hbox{ if } k<2(t-m) \\
    H_{k-2(t-m)}(X; \Q), & \hbox{ if } k\geq 2(t-m) \\
\end{array}%
\right.$$ Therefore, if $t\geq m$, since all maps in the chain of
isomorphisms
$$L^tH^k(X;\Q)\cong \pi_{2t-k}Z_0(X\times \C^{t-m})\otimes \Q\cong H^{BM}_{2t-k}(X\times
\C^{t-m}; \Q)\cong H_{2m-k}(X; \Q)\cong H^k(X; \Q)$$ are natural,
their composite is $\Phi^{t, k}$. And from the commutative diagram
of the $S$-map, we see that the $S$-map is also an isomorphism.
\end{proof}

Suppose that $(\ , \ )$ is a symmetric bilinear form on a finite
dimensional vector space $V$ over $\Q$. The signature of $(\ , \ )$
is the number of positive eigenvalues minus the number of negative
eigenvalues in a matrix representation of $(\ , \ )$.

From the natural transformation $\Phi^{t, k}\otimes \F: L^tH^k(X;
\F)\longrightarrow H^k(X; \F)$, we define
$$\widetilde{L}^tH^k(X; \F)=Im(\Phi^{t, k}\otimes \F)\subset H^k(X; \F).$$

\begin{definition}(morphic signatures)
Suppose the dimension of $X$ is $m=2n$. For $t\geq n$, we define the
$t$-th morphic signature of $X$, denoted by $\sigma_t$, to be the
signature of the symmetric bilinear form:
$$(\ , \ ): \widetilde{L}^tH^m(X; \Q)\otimes \widetilde{L}^tH^m(X; \Q)\longrightarrow
\widetilde{L}^{2t}H^{2m}(X; \Q)=\Q.$$
\end{definition}

For $t=m$, since $\widetilde{L}^{m}H^{m}(X; \Q)=H^{m}(X; \Q)$ and
the cup product in morphic cohomology in this case is just the usual
cup product of singular cohomology, $\sigma_{m}$ is the usual
signature of $X$. So we have a sequence of signatures $\sigma_{m},
\sigma_{m-1}, ..., \sigma_n$ which reveals more and more algebraic
information of $X$.

\section{The Morphic Conjectures}
Let $a, b$ be two nonnegative integers. Define
$$EH^a(X;
\F)=\widetilde{L}^aH^0(X; \F)\oplus\widetilde{L}^{a+1}H^2(X;
\F)\oplus \cdots \oplus \widetilde{L}^{a+m}H^{2m}(X; \F),$$
$$OH^b(X; \F)=\widetilde{L}^bH^1(X; \F)\oplus
\widetilde{L}^{b+1}H^3(X; \F)\oplus \cdots \oplus
\widetilde{L}^{b+m-1}H^{2m-1}(X; \F)$$ where $E$ and $O$ stand for
even and odd respectively. In particular, $EH^0(X; \Q)$ is the ring
of rational algebraic cohomology classes on $X$. Define
$$LH^{a, b}(X; \F)=EH^a(X; \F)\oplus OH^b(X; \F).$$

Let $\Omega\in L^1H^2(X)$ be a class coming from a hyperplane
section on $X$. Define an operation
$$\mathcal{L}: L^tH^k(X) \longrightarrow L^{t+1}H^{k+2}(X)$$
by $\mathcal{L}(\alpha)=\Omega\cdot \alpha$. Under the
transformation $\Phi^{*, *}$, $\mathcal{L}$ carries over to the
standard Lefschetz operator $L$. The operators induced by
$\mathcal{L}$ on $EH^a(X; \F), OH^b(X; \F)$ and $LH^{a, b}(X; \F)$
are simply the restriction of $L$ to these spaces, and these spaces
are $L$-invariant. By abuse of notation, we use $\mathcal{L}$ to
denote the restriction of $L$ to these spaces.

Recall that there is a standard Hermitian inner product on
$\mathscr{A}^{p, q}(X)$, the $(p, q)$-forms on $X$, called the
Hodge inner product defined by
$$<\alpha, \beta>=\int_X\alpha \wedge *\bar{\beta}$$
where $*$ is the Hodge star operator. Let $\Lambda$ be the adjoint
of $L$ with respect to the Hodge inner product. Since $L, \Lambda$
commute with the Laplacian, they can be defined on the harmonic
spaces. From the Hodge theorem we know that the $(p, q)$-cohomology
group of $X$ is isomorphic to the space of $(p, q)$-harmonic forms.
The Hodge inner product induces a Hermitian inner product in
harmonic spaces which we also call the Hodge inner product. Restrict
the Hodge inner product to $EH^a(X; \F), OH^b(X; \F)$ and $LH^{a,
b}(X; \F)$ respectively and let $\lambda$ be the adjoint of
$\mathcal{L}$ with respect to the Hodge inner product.

\begin{conjecture}(morphic conjectures)
The morphic conjecture on $EH^a(X; \F), OH^b(X; \F)$ and $LH^{a,
b}(X; \F)$ respectively is the assertion that $\lambda$ is the
restriction of $\Lambda$ on them respectively.
\end{conjecture}

It is not difficult to see that if a morphic conjecture holds for
$\Q$-coefficients, it also holds for $\R$- and $\C$-coefficients. So
most of the time we will only work with $\Q$-coefficients.

\begin{definition}
Consider the cup product pairing:
$$L^tH^k(X; \Q)\otimes L^{t+m-k}H^{2m-k}(X; \Q) \longrightarrow
L^{2t+m-k}H^{2m}(X; \Q)$$ For $\alpha \in L^tH^k(X; \Q)$, we say
that $\alpha$ is morphic numerically equivalent to $0$ if $\alpha
\wedge \beta=0$ for all $\beta \in L^{t+m-k}H^{2m-k}(X; \Q)$. The
class $\alpha$ is said to be morphic homologically equivalent to $0$
if $\Phi(\alpha)=0$ where $\Phi: L^tH^k(X; \Q)\longrightarrow H^k(X;
\Q)$ is the natural transformation. We use $MNE$ for morphic
numerical equivalence and $MHE$ for morphic homological equivalence.
\end{definition}

Let $Alg_k(X; \Q)$ be the group of $k$-cycles with rational
coefficients on $X$ quotient by algebraic equivalence and let
$Alg^k(X; \Q)=Alg_{m-k}(X; \Q)$. We recall that a class $\alpha\in
Alg^k(X; \Q)$ is said to be numerically equivalent to zero if
$\alpha \bullet \beta=0$ for all $\beta \in Alg^{m-k}(X; \Q)$ where
$\bullet$ is the intersection product and $\alpha$ is said to be
homologically equivalent to zero if under the cycle map $\gamma:
Alg^k(X; \Q) \longrightarrow H^{2k}(X; \Q)$, $\alpha$ is sent to
zero (see e.g \cite{Kle1}). By the Friedlander-Lawson duality
theorem, we can identify $L^tH^{2t}(X; \Q)$ with $Alg^t(X; \Q)$ for
$0\leq t \leq m$ (see \cite[Theorem 5.1]{FL3}).

From what we explain before, we have the following result.
\begin{proposition}
\begin{enumerate}
\item On $\oplus^m_{t=0}L^tH^{2t}(X;\Q)$ morphic numerical equivalence is
same as numerical equivalence and morphic homological equivalence is
same as homological equivalence.

\item If $\alpha$ is morphic homologically equivalent to zero, then
$\alpha$ is morphic numerically equivalent to zero.
\end{enumerate}
\end{proposition}

\begin{proposition}
\begin{enumerate}
\item $dim\widetilde{L}^{a+t}H^{2t}(X; \Q) \leq
dim\widetilde{L}^{a+m-t}H^{2m-2t}(X; \Q)$ for $t\leq \lfloor
\frac{m}{2} \rfloor$

\item $dim\widetilde{L}^{b+t}H^{2t-1}(X; \Q) \leq
dim\widetilde{L}^{b+m-t+1}H^{2m-2t+1}(X; \Q)$ for $t\leq \lfloor
\frac{m}{2} \rfloor$
\end{enumerate}
\end{proposition}

\begin{proof}
We have the following commutative diagrams:
$$\xymatrix{ L^{a+t}H^{2t} \ar[rr]^-{\mathcal{L}^{m-2t}} \ar[d]& &
L^{a+m-t}H^{2m-2t} \ar[d] & & L^{b+t}H^{2t-1}
\ar[rr]^-{\mathcal{L}^{m-2t+1}} \ar[d] & & L^{b+m-t+1}H^{2m-2t+1} \ar[d]\\
H^{2t} \ar[rr]^-{L^{m-2t}} & & H^{2m-2t} & & H^{2t-1}
\ar[rr]^-{L^{m-2t+1}} & & H^{2m-2t+1}}$$

By the Hard Lefschetz theorem, $L^{m-2t}$ and $L^{m-2t+1}$ are
isomorphisms for $t\leq \lfloor \frac{m}{2} \rfloor$ , so we have
the conclusions.
\end{proof}

Let $\mathscr{A}$ be any one of $EH^a(X; \Q), OH^b(X; \Q)$ and
$LH^{a, b}(X; \Q)$. Let $\widetilde{\mathscr{A}}$ be the direct sum
of all morphic cohomology group $L^tH^k(X; \Q)$ such that
$\widetilde{L}^tH^k(X; \Q)$ is a direct summand of $\mathscr{A}$.

\begin{proposition}
The followings are equivalent:
\begin{enumerate}
\item MNE=MHE on $\widetilde{\mathscr{A}}$.

\item $dim\widetilde{L}^tH^k(X;
\Q)=dim\widetilde{L}^{t+m-k}H^{2m-k}(X; \Q)$ for
$\widetilde{L}^tH^k(X; \Q) \subset \mathscr{A}$.

\item If $\alpha \in \widetilde{L}^tH^k(X; \Q) \subset \mathscr{A} $ for $k\leq m$ and
$\alpha=\sum_{r\geq 0} L^r\alpha_r$ is the Lefschetz decomposition
of $\alpha$, then $\alpha_r\in \widetilde{L}^{t-r}H^{k-2r}(X;
\Q)\subset \mathscr{A}$, for $r\geq 0$.

\item If $\alpha \in \mathscr{A}$ then $*\alpha \in \mathscr{A}$.

\item If $\alpha \in \mathscr{A}$ then $\Lambda \alpha \in \mathscr{A}$.

\item The morphic conjecture holds on $\mathscr{A}$.
\end{enumerate}
\label{equivalent forms}
\end{proposition}

\begin{proof}
We are going to show that $1\rightarrow 2 \rightarrow 3
\rightarrow 4 \rightarrow 5 \rightarrow 6 \rightarrow 2$ and $4
\rightarrow 1$.

$1 \rightarrow 2$: We consider only the case $\mathscr{A}=EH^a(X;
\Q)$ since similar argument applies for $\mathscr{A}=OH^b(X; \Q)$
and thus for $\mathscr{A}=LH^{a, b}(X; \Q)$. For $t\leq \lfloor
\frac{m}{2} \rfloor$, consider the commutative diagram
$$\xymatrix{ L^{a+t}H^{2t}(X; \Q)\otimes L^{a+m-t}H^{2m-2t}(X; \Q)\ar[d] \ar[r] &
L^{2a+m}H^{2m}(X; \Q) \ar[d]\\
H^{2t}(X; \Q)\otimes H^{2m-2t}(X; \Q) \ar[r] & H^{2m}(X; \Q)}$$

If $\widetilde{L}^{a+m-t}H^{2m-2t}(X; \Q)=0$, by the proposition
above, $\widetilde{L}^{a+t}H^{2t}(X; \Q)=0$. So we may assume that
$\widetilde{L}^{a+m-t}H^{2m-2t}(X; \Q)\neq 0$. Let $\alpha \in
L^{a+m-t}H^{2m-2t}(X; \Q)$ such that $\Phi(\alpha)\neq 0$. If
MNE=MHE on $\widetilde{\mathscr{A}}$, then there is a $\beta\in
L^{a+t}H^{2t}(X; \Q)$ such that $(\beta, \alpha)\neq 0$ where $(\ ,
\ )$ is the cup product pairing in morphic cohomology. Thus
$(\Phi(\beta), \Phi(\alpha))\neq 0$. Consequently, the cup product
pairing $(\ , \ )$ in singular cohomology restricted to
$\widetilde{L}^{a+t}H^{2t}(X; \Q)\otimes
\widetilde{L}^{a+m-t}H^{2m-2t}(X; \Q)$ is nondegenerate. It follows
that $dim\widetilde{L}^{a+t}H^{2t}(X; \Q) =
dim\widetilde{L}^{a+m-t}H^{2m-2t}(X; \Q)$.

$2 \rightarrow 3$: From the commutative diagram
$$\xymatrix{L^{t-1}H^{k-2}(X; \Q) \ar[d] \ar[rr]^-{\mathcal{L}^{m-k+2}} & &
L^{t+m-k+1}H^{2m-k+2}(X; \Q) \ar[d]\\
H^{k-2}(X; \Q) \ar[rr]^{L^{m-k+2}} & & H^{2m-k+2}(X; \Q)}$$ we see
that $L^{m-k+2}$ maps $\widetilde{L}^{t-1}H^{k-2}(X; \Q)$
injectively into $\widetilde{L}^{t+m-k+1}H^{2m-k+2}(X; \Q)$. The
assumption $dim \widetilde{L}^{t-1}H^{k-2}(X;
\Q)=dim\widetilde{L}^{t+m-k+1}H^{2m-k+2}(X; \Q)$ implies that
$L^{m-k+2}$ restricted to $\widetilde{L}^{t-1}H^{k-2}(X; \Q)$ is
an isomorphism. Let $\alpha=\sum_{r\geq 0}L^r\alpha_r\in
\widetilde{L}^tH^k(X; \Q)$ be the Lefschetz decomposition of
$\alpha$. We prove by induction on the length of the Lefschetz
decomposition. Since $L^{m-k+2}(\sum_{r\geq
1}L^{r-1}\alpha_r)=L^{m-k+1}(\alpha) \in
\widetilde{L}^{t+m-k+1}H^{2m-k+2}(X; \Q)$, we have $\sum_{r\geq
1}L^{r-1}\alpha_r\in \widetilde{L}^{t-1}H^{k-2}(X; \Q)$. By
induction hypothesis, $\alpha_r\in \widetilde{L}^{t-r}H^{k-2r}(X;
\Q)$ for $r\geq 1$. But $\alpha_0=\alpha-L(\sum_{r\geq
1}L^{r-1}\alpha_r)\in \widetilde{L}^tH^k(X; \Q)$. This completes
the proof.

$3 \rightarrow 4$: Suppose that $\alpha \in \widetilde{L}^tH^k(X;
\Q)$ and $\alpha=\sum_{r\geq 0}L^r\alpha_r$ is the Lefschetz
decomposition of $\alpha$. By some calculation we get
$*L^j\beta=(-1)^{\frac{k(k+1)}{2}}\frac{j!}{(m-k-j)!}L^{m-k-j}\beta$
for $\beta \in B^k$ as a formula in Definition \ref{star}. From the
assumption $\alpha_r\in \mathscr{A}$, we have $L^{m-k+r}\alpha_r\in
\mathscr{A}$, for $r\geq 0$. Thus $*\alpha \in \mathscr{A}$.

$4 \rightarrow 5$: From the formula $\Lambda=*L*$ as in Definition
\ref{star}, we have the conclusion immediately.

$5 \rightarrow 6$: Since $\lambda=\pi\circ \Lambda$ where
$\pi:\oplus^{2m}_{k=0}H^k(X; \Q)\longrightarrow \mathscr{A}$ is the
projection, from the assumption $\pi \circ
\Lambda|_{\mathscr{A}}=\Lambda|_{\mathscr{A}}$, we have
$\lambda=\Lambda|_{\mathscr{A}}$. Therefore, the morphic conjecture
holds on $\mathscr{A}$.

$6 \rightarrow 2$: By the Hard Lefschetz theorem, $\Lambda^{m-k}:
H^{2m-k}(X; \Q)\longrightarrow H^k(X; \Q)$ is an isomorphism for
$k\leq m$. Therefore if $\lambda=\Lambda|_{\mathscr{A}}$,
$\Lambda^{m-k}(\widetilde{L}^{t+m-k}H^{2m-k}(X; \Q))\subset
\widetilde{L}^tH^k(X; \Q)$ which implies that they have the same
dimension.

$4 \rightarrow 1$: Suppose that $\alpha\in L^tH^k(X; \Q)$ is
morphic numerically equivalent to zero. If $\Phi(\alpha)\neq 0$
then
$$(\Phi(\alpha), *\Phi(\alpha))=\int_X \Phi(\alpha)\wedge
*\Phi(\alpha)=<\Phi(\alpha), \Phi(\alpha)>\neq 0$$ But by the
hypothesis $*\Phi(\alpha)\in \widetilde{L}^{t+m-k}H^{2m-k}(X;
\Q)$, so we can find $\beta\in L^{t+m-k}H^{2m-k}(X; \Q)$ such that
$\Phi(\beta)=*\Phi(\alpha)$. Then $(\alpha, \beta)=(\Phi(\alpha),
\Phi(\beta))\neq 0$ which contradicts to the assumption. Thus
$\alpha$ is morphic homologically equivalent to zero.
\end{proof}

In particular, GSCA is equivalent to the morphic conjecture on
$EH^0(X)$.

\begin{definition}
For $\beta \in L^rH^{2r}(X\times X; \Q)$, $\beta$ induces a map
$$\beta_*: L^tH^k(X; \Q) \longrightarrow L^{t-(m-r)}H^{k-2(m-r)}(X; \Q)$$ defined by
$$\beta_*(\alpha)=\mathscr{D}^{-1}(q_*(\mathscr{D}(p^*\alpha)\bullet \mathscr{D}(\beta)))$$
where $p, q:X\times X\longrightarrow X$ are the projections to the
first and second factor respectively, $\mathscr{D}$ is the
Friedlander-Lawson duality map and $\bullet$ is the intersection
product in Lawson homology. An endomorphism $f: \mathscr{A}
\longrightarrow \mathscr{A}$ is said to be algebraic if there is a
$\beta \in \oplus^{m}_{r=0}L^rH^{2r}(X\times X; \Q)$ such that
$\beta_*=f$.
\end{definition}

We note that this definition is equivalent to the definition in
\cite{Groth}.

\begin{proposition}
If the Grothendieck standard conjecture B holds on $X$, then it
implies all the equivalent statements in Proposition \ref{equivalent
forms} for $\mathscr{A}$ equals to any of $EH^a(X; \Q), OH^b(X; \Q)$
or $LH^{a, b}(X; \Q)$.
\end{proposition}

\begin{proof}
If the GSCB holds on $X$, then $\Lambda$ is an algebraic operator,
thus there exists a cycle $\beta\in L^{m-1}H^{2(m-1)}(X\times X;
\Q)$ such that $\Lambda=\beta_*$. For $L^tH^k(X; \Q)$ a direct
summand of $\widetilde{\mathscr{A}}$, $\beta_*(L^tH^k(X; \Q))
\subset L^{t-1}H^{k-2}(X; \Q)$, thus $\Lambda$ is an endomorphism of
$\mathscr{A}$. By the fifth statement in Proposition \ref{equivalent
forms}, the morphic conjecture holds on $\mathscr{A}$.
\end{proof}

Hence all the morphic conjectures are true for abelian varieties,
varieties of complete intersection and Grassmannians.

By assuming the GSCB, we answer a question of Friedlander and Lawson
in rational coefficients (see \cite{FL1}, Question 9.7).

\begin{theorem}
If the Grothendieck standard conjecture B holds on $X$, then the map
$$\Phi^{t, k}:L^tH^k(X; \Q) \longrightarrow H^k(X; \Q)$$ is
surjective whenever $t\geq k$. \label{surjectivity}
\end{theorem}

\begin{proof}
By Proposition \ref{greater than the dimension}, it is true if
$t\geq m$. So we assume that $t<m$. Then $k<m$. If the GSCB holds on
$X$, we have all the morphic conjectures. Thus from Proposition
\ref{equivalent forms}, the dimension of
$\widetilde{L}^{m+t-k}H^{2m-k}(X; \Q)$ is same as the dimension of
$\widetilde{L}^tH^k(X; \Q)$. Since $m+t-k\geq m$,
$\widetilde{L}^{m+t-k}H^{2m-k}(X; \Q)=H^{2m-k}(X; \Q)$, and by the
Hard Lefschetz theorem, we have $\widetilde{L}^tH^k(X; \Q)=H^k(X;
\Q)$. Therefore $\Phi^{t, k}$ is surjective.
\end{proof}

\section{Topological Filtration And Arithmetic Filtration}
In \cite{FM}, Friedlander and Mazur defines two filtrations on the
singular homology groups of a projective variety. One is formed by
taking the images of the natural transformations from Lawson
homology to singular homology and the other one is the homological
version of the Grothendieck's arithmetic filtration. The first
filtration is called the topological filtration (denoted by
$T_rH_n(X; \Q)$) and the second one is called the geometric
filtration (denoted by $G_rH_n(X; \Q)$). Friedlander and Mazur
conjecture that these two filtrations are equal.

\begin{conjecture}(Friedlander-Mazur)
Let $j, n$ be non-negative integers. For any smooth projective
variety $X$,
$$T_jH_n(X; \Q)=G_jH_n(X; \Q)$$
\end{conjecture}

In the following, we define a filtration from morphic cohomology
and reformulate the Friedlander-Mazur conjecture as an equality
between this filtration and the Grothendieck's arithmetic
filtration.

For a variety $Y$, let $\gamma:\widetilde{Y}\longrightarrow Y$ be
a desingularization of $Y$. Recall that the arithmetic filtration
(coniveau filtration) $\{N^pH^*(X; \Q)\}_{p\geq 0}$ of $H^*(X;
\Q)$ is given by
$$N^pH^l(X; \Q)=\{\mbox{ Gysin images } \gamma_*: H^{l-2q}(\widetilde{Y};
\Q) \longrightarrow H^l(X; \Q)|Y \subset X, \mbox{ codim }_XY=q
\mbox{ (pure) }, q\geq p \}$$ (see \cite{Lew}, page 87 for details);
and recall that the niveau filtration $\{N_pH_*(X; \Q)\}_{p\geq 0}$
of $H_*(X; \Q)$ is defined by
$$N_pH_i(X; \Q)=\mbox{ span }\{\mbox{ images } i_*:H_i(Y; \Q)
\longrightarrow H_i(X; \Q)|i: Y \hookrightarrow X, \mbox{ dim Y}
\leq p\}$$

Define the topological filtration $\{T^pH^*(X; \Q)\}$ to be
$$T^pH^l(X; \Q)=\{ \mbox{ images } \Phi^{p, l}: L^pH^1(X; \Q)
\longrightarrow H^l(X; \Q)\}$$ where $\Phi^{p, l}$ is the natural
transformation from morphic cohomology to singular cohomology.

If $X$ is a smooth projective variety of dimension $m$, it is not
difficult to see that $G_rH_n(X; \Q)=N_{n-r}H_n(X; \Q)$ and
$N_pH_i(X; \Q)\cong N^{m-p}H^{2m-i}(X; \Q)$ by the Poincar\'e
duality. By \cite[Theorem 5.9]{FL3}, $T^rH^n(X; \Q) \cong
T_{m-r}H_{2m-n}(X; \Q)$. It is proved in \cite[7.5, Corollary
3]{FM}, that $T_rH_n(X; \Q) \subset G_rH_n(X; \Q)$. The
cohomological version of this result is the containment
$T^{l-p}H^l(X; \Q) \subset N^pH^l(X; \Q)$ and the cohomological
version of the Friedlander-Mazur conjecture is the following
conjecture.

\begin{conjecture}
For nonnegative integers $l, p$, $T^{l-p}H^l(X; \Q)=N^pH^l(X;
\Q)$.
\end{conjecture}

Recall that the generalized Hodge conjecture is the assertion that
$N^pH^l(X; \Q)=F^p_hH^l(X; \Q)$ for all $p, l$ where
$\{F^P_hH^*(X; \Q)\}$ is the rational Hodge filtration (see
\cite{Lew}, page 88). If the Friedlander-Mazur conjecture holds,
ideally, it would give a more concrete picture about the
arithmetic filtration.

As a consequence of Theorem \ref{surjectivity}, we have some
evidence for the Friedlander-Mazur conjecture.

\begin{corollary}
If the Grothendieck standard conjecture B holds on a smooth
projective variety $X$. Then
$$T^tH^k(X; \Q)=N^0H^k(X; \Q)=F^0_hH^k(X; \Q)=H^k(X; \Q)$$ for
$t\geq k$.
\end{corollary}

In the following, we give a simple proof of the Friedlander's result
in \cite[Proposition 4.2]{F3}.

\begin{proposition}
Suppose that $X$ is a smooth projective variety. If the
Grothendieck standard conjecture B is valid for a resolution of
singularities of each irreducible subvariety $Y\subset X$ of
codimension $\geq p$. Then
$$N^pH^l(X; \Q)=T^{l-p}H^l(X; \Q)$$
\end{proposition}

\begin{proof}
Suppose that $Y\subset X$ is a subvariety of codimension $p'\geq
p$. Let $\sigma:\widetilde{Y} \longrightarrow Y$ be a
desingularization and the GSCB holds on $\widetilde{Y}$. Consider
the following commutative diagram:
$$\xymatrix{L^{l-p-p'}H^{l-2p'}(\widetilde{Y}; \Q)
\ar[r]^-{\sigma_*} \ar[d]_{\Phi^{l-p-p', l-2p'}} & L^{l-p}H^l(X; \Q) \ar[d]^{\Phi^{l-p, l}} \\
H^{l-2p'}(\widetilde{Y}; \Q) \ar[r]^-{\sigma_*} & H^l(X; \Q)}$$ By
Theorem \ref{surjectivity}, $\Phi^{l-p-p', l-2p'}$ is surjective.
Therefore the image of $\sigma_*$ is contained in the image of
$\Phi^{l-p, l}$. Therefore, $N^pH^l(X; \Q) \subset T^{l-p}H^l(X)$.
\end{proof}

Since the GSCB holds for smooth projective varieties of dimension
$\leq 2$, we have the following result.

\begin{corollary}
Suppose that $X$ is a smooth projective variety of dimension less
than or equal to 3. Then $N^pH^l(X)=T^{l-p}H^l(X)$ for all $p, l$.
\end{corollary}

\section{Abstract Hodge index theorem}
\begin{definition}
Let $V=\oplus^{2m}_{k=0}H^k$ where each $H^k$ is a finite
dimensional vector space over $\Q$ and let $V_{\F}=V\otimes \F$,
$H^k_{\F}=H^k\otimes \F$. Let $h=\sum^{2m}_{k=0}(k-m)\pi_k$ where
$\pi_k:V_{\C} \rightarrow H^k_{\C}$ is the projection. $V$ is called
a Lefschetz algebra if
\begin{enumerate}
\item there is an inner product $< , >: V_{\R}\times V_{\R} \rightarrow \R$
which induces a hermitian inner product $< , >: V_{\C}\times V_{\C}
\rightarrow \C$ defined by $<a\otimes \mu, b\otimes \lambda>:=\mu
\overline{\lambda}<a, b>$.

\item There is an endomorphism $L:V \rightarrow V$ of degree 2 with
adjoint $\Lambda$ such that $L, \Lambda$ and $h$ define a
$sl_2(\C)$-action on $V_{\C}$ in the following way:
$$[\Lambda, L]=h, [h, \Lambda]=2\Lambda, [h, L]=-2L$$
\end{enumerate}
\end{definition}

\begin{definition}
Let $B^k:=ker\Lambda:H^k_{\C} \rightarrow H^{k-1}_{\C}$ be the
primitive space. For $\alpha\in B^k$, define
$$\overline{*}L^j\alpha:=(-1)^{\frac{k(k+1)}{2}}\frac{j!}{(m-k-j)!}L^{m-k-j}\overline{\alpha}$$
and define
$$\Lambda L^j\alpha:=j(m-k-j+1)L^{j-1}\alpha$$
\label{star}
\end{definition}

\begin{proposition}
For a Lefschetz algebra $V$ as above, we have the following
properties:
\begin{enumerate}
\item There is a Strong Lefschetz theorem:
$$L^{m-k}: H^k \overset{\cong}{\longrightarrow} H^{2m-k}$$

\item There is a Lefschetz decomposition:
for $a\in H^k$, $$a=\sum_{j\geq max(0, k-m)}  L^j\alpha_j$$ where
$\alpha_j\in B^{k-2j}$.

\item The Lefschetz decomposition is orthogonal with respect to $<
, >$.

\item $\overline{*}^2=id$, $\overline{*}$ is conjugate self-adjoint,
i.e., $<\alpha, \overline{*}\beta>=\overline{<\overline{*}\alpha,
\beta>}$.

\item $\Lambda=\overline{*}L\overline{*}$.
\end{enumerate}
\label{orthogonal}
\end{proposition}

\begin{proof}
(1) and (2) follow from the properties of $sl_2(\C)$-action (see
\cite[Theorem 11.15]{Lew}). Let $L^k\alpha \in L^kB^{n-2k}, L^s\beta
\in L^sB^{n'-2s}$ and $k\geq s$. By the relation $[\Lambda, L]=h$,
we have $<L\alpha, L\beta>=<\Lambda L \alpha, \beta>=<h\alpha +
L\Lambda \alpha, \beta>=<h\alpha, \beta>=c<\alpha, \beta>$ where $c$
is a constant. Hence if $k>s$, $<L^k\alpha, L^s
\beta>=c<L^{k-s}\alpha, \beta>=c<L^{k-s-1}\alpha, \Lambda \beta>=0$.
Hence the decomposition is orthogonal with respect to $< , >$. (4)
and (5) follow from some simple calculations.
\end{proof}

\begin{definition}
Let $V=\oplus^{2n}_{t=0}H^{2t}$ be a Lefschetz algebra. Suppose that
$V$ is endowed with the following structures:

\begin{enumerate}
\item Each $H^{2t}_{\C}=\oplus_{p+q=2t}H^{p, q}$
has a Hodge structure of weight $2t$ such that the decomposition is
orthogonal with respect to $< , >$.

\item The Hodge structure is compatible with the $sl_2(\C)$-action,
i.e.,
$$L^k:H^{p, q} \rightarrow H^{p+k, q+k}$$
for any $p, q, k$.

\item Let $B^{p, q}:=ker \Lambda: H^{p, q} \rightarrow H^{p-1,
q-1}$ and define
$$Q(\alpha, \beta):=<L^{n-r}\alpha, L^{n-r}\overline{\beta}>$$
for $\alpha\in B^{p, q}, \beta \in B^{q, p}$ and $2r=p+q$. There are
the Hodge-Riemann bilinear relations:
\begin{enumerate}
\item $Q(B^{p, q}, B^{s, t})=0$ if $s\neq q$.
\item $(-1)^{r+q}Q(\xi, \overline{\xi})>0$ if $0\neq \xi\in B^{p,
q}$ where $2r=p+q$.
\end{enumerate}
\end{enumerate}
Then we say that $V$ is a Hodge-Lefschetz algebra.
\end{definition}

From now on, our $V$ denote a Hodge-Lefschetz algebra as above.
Since the $sl_2(\C)$-action is compatible with the Hodge structure,
it reduces to an $sl_2(\C)$-action on $V^{a, b}_{\C}=\oplus^{min(a,
b)}_{k=-min(a, b)} H^{a+k, b+k}$. Hence for all $p, q$, we have a
Lefschetz decomposition
$$H^{p, q}=B^{p, q}\oplus LB^{p-1, q-1} \oplus L^2B^{p-2, q-2}\oplus
\cdots \oplus L^rB^{p-r, q-r}$$ where $r=min(p, q)$. Similar to the
proof of Proposition \ref{orthogonal}, we have the orthogonality of
the decomposition.

\begin{proposition}
The decomposition
$$H^{2t}_{\C}=\bigoplus_{p+q=2t}\bigoplus_{0\leq k \leq min(p, q)} L^kB^{p-k,
q-k}$$ is orthogonal with respect to $< , >$.
\end{proposition}

Let $h^{p, q}=dim_{\C}H^{p, q}_{\C}$. We follow \cite[Theorem
15.8.2]{Hirz} to give a proof of the Hodge index theorem.
\begin{theorem}(Abstract Hodge index theorem)
Suppose that $V$ is a Hodge-Lefschetz algebra as above. Define
$(\alpha, \beta):=<\alpha, \overline{*}\beta>$ on $H^{2n}$. Then the
signature $\sigma$ of $(\ , \ )$ is $\sum_{p, q}(-1)^qh^{p, q}$.
\label{Abstract Hodge index theorem}
\end{theorem}

\begin{proof}
\begin{enumerate}
\item Since $\overline{*}_{|_V}$ is self-adjoint and $< ,  >_{|_{V\times V}}$ is symmetric, $(\ ,\ )$ is a
symmetric bilinear form.

\item Let $E^{p, q}_k$ be the vector space consisting of
$L^k(a+\overline{a})$ for $a \in B^{p-k, q-k}$. Then $E^{p, q}_k$ is
a real vector space. We have
$$H^{2n}_{\R}=\bigoplus_{p+q=2n}\bigoplus_{0\leq k \leq min(p, q)}E^{p,
q}_k$$

\item The decomposition above is orthogonal with respect to the
Hodge inner product and the quadratic form $(-1)^{q+k}(\ ,\ )$ is
positive definite when restricted to $E^{p, q}_k$.

\begin{proof}
For $\alpha\in E^{p, q}_k$ where $p+q=2n$ and $p\neq q$, let
$\alpha=L^ka$ where $a=b+\overline{b}$ is real, $b\in B^{p-k, q-k}$,
by a simple calculation, we have
$\overline{*}L^kb=(-1)^{n-k}L^k\overline{b}$ then $(\alpha,
\alpha)=2<L^kb, \overline{*}L^kb>=2(-1)^{n-k}<L^kb,
L^kb>=2(-1)^{n-k}Q(b, \overline{b})$. Hence by the Hodge-Riemann
bilinear relation, $(-1)^{q+k}(b, b)=2(-1)^{(n-k)+(q-k)}Q(b,
\overline{b})>0$ if $b\neq 0$. Similarly, if $p=q$, let
$\alpha=L^kb$ where $b\in B^{n-k, n-k}$ and $b=\overline{b}$. Then
$(\alpha, \alpha)=<\alpha, *\alpha>=(-1)^{n+k}<L^kb,
L^kb>=(-1)^{n+k}Q(b, \overline{b})$. Hence, $(-1)^{n+k}(\alpha,
\alpha)=Q(b,\overline{b})>0$ if $b\neq 0$.
\end{proof}

\item Therefore, $$\sigma=\sum_{\substack{p+q=2n \\ k\leq p \leq q}}(-1)^{q+k}dim_{\R}E^{p,
q}_k$$

\item $$\sigma=\sum_{\substack{p+q=2n\\k\leq min(p,
q)}}(-1)^{q+k}dim_{\C}L^kB^{p-k, q-k}$$

\begin{proof}
The real dimension of $E^{p, q}_k$ is $dim_{\C}L^kB^{p-k,
q-k}+dim_{\C}L^kB^{p-k, q-k}$ for $p<q$ and $dim_{\R}E^{n,
n}_k=dim_{\C}L^kB^{n-k, n-k}$.
\end{proof}

\item
$h^{p-k, q-k}-h^{p-k-1, q-k-1}=dim_{\C} B^{p-k, q-k}=dim_{\C}
L^kB^{p-k, q-k}$ for $p+q\leq 2n$.

\item For $p+q=2n$, by the Hard Lefschetz theorem, we have $h^{p-k-1, q-k-1}=h^{p+k+1, q+k+1}$ and from the Hodge
structures, $h^{r, k}=h^{k, r}=h^{2n-r, 2n-k}$.

\item $$\sigma=\sum_{\substack{k\geq 0\\p+q=2n}}(-1)^{q-k}h^{p-k,
q-k}+\sum_{\substack{k\geq 0\\p+q=2n}}(-1)^{q+k+1}h^{p+k+1, q+k+1}$$
$$=\sum_{p+q\leq 2n}(-1)^qh^{p, q}+\sum_{p+q>2n}(-1)^qh^{p, q}
=\sum_{p, q}(-1)^qh^{p, q}$$
\end{enumerate}
\end{proof}

\section{Morphic Hodge Index Theorem}
Let us use $H^*(X; \C)$ to denote the cohomology ring of $X$. Let
$h=\sum^{2m}_{k=0}(m-k)Pr_k$ where $Pr_k:H^*(X; \C) \longrightarrow
H^k(X; \C)$ projects a form to its $k$-component. The
$sl_2(\C)$-structure on $H^*(X; \C)$ is given by
$$[\Lambda, L]=h, [h, \Lambda]=2\Lambda, [h, L]=-2L$$

Let $\mathscr{A}$ be any of $EH^a(X; \C), OH^b(X; \C)$ or $LH^{a,
b}(X;
\C)$ and $$\gamma(\mathscr{A}, p, q)=\left\{%
\begin{array}{ll}
    a+\frac{p+q}{2}, & \hbox{ if } p+q \mbox{ is even } \\
    b+\frac{p+q-1}{2}, & \hbox{ if } p+q \mbox{ is odd }.\\
\end{array}%
\right.
$$
If it is clear from the context what $\mathscr{A}$ is, to simplify
our notation, we will just write $\gamma(p, q)$ for
$\gamma(\mathscr{A}, p, q)$. Let $\mathcal{L}$ be the restriction of
$L$ to $\mathscr{A}$ and $\lambda$ be the adjoint of $\mathcal{L}$
with respect to the Hodge inner product restricted to $\mathscr{A}$.

\begin{proposition}
Assume that the morphic conjecture holds on $\mathscr{A}$, then
$\mathscr{A}$ is a $sl_2(\C)$-submodule of $H^*(X; \C)$ thus
\begin{enumerate}
\item $\mathscr{A}$ has a sub-Lefschetz decomposition, i.e.,
if $\widetilde{L}^tH^k(X; \C)$ is a direct summand of
$\mathscr{A}$, then
$$\widetilde{L}^tH^k(X; \C)=\oplus_{r\geq \mbox{max}\{0,
k-m\}}\mathcal{L}^rB^{m-2r}$$ where
$B^k=Ker\mathcal{L}^{m-k+1}:\widetilde{L}^tH^k(X; \C)
\longrightarrow \widetilde{L}^{t+m-k+1}H^{m-k+2}(X; \C)$ is the
primitive group. Furthermore, this decomposition is compatible
with the sub-Hodge structure, i.e., if $B^{p,
q}=KerL^{m+1-p-q}:H^{p, q}_{\gamma(p, q)}(X)\longrightarrow
H^{m+1-q, m+1-p}_{\gamma(m+1-q, m+1-p)}(X)$, then
$$H^{p, q}_{\gamma(p, q)}(X)=\oplus_{r\geq \mbox{max}\{0,
k-m\}}\mathcal{L}^rB^{p-r, q-r}$$

\item $B^k=ker\lambda:\widetilde{L}^tH^k(X; \C) \longrightarrow
\widetilde{L}^{t-1}H^{k-2}(X; \C)$ and $B^{p, q}=ker\lambda: H^{p,
q}_t(X) \longrightarrow H^{p-1, q-1}_{t-1}(X)$ where
$\widetilde{L}^tH^k(X; \C)$ is a direct summand of $\mathscr{A}$.

\item We have the Hard Lefschetz theorem, i.e.,
$$\mathcal{L}^k: \widetilde{L}^tH^{m-k}(X; \C) \longrightarrow
\widetilde{L}^{t+k}H^{m+k}(X; \C)$$ is an isomorphism where
$\widetilde{L}^tH^{m-k}(X; \C)$ is a direct summand of
$\mathscr{A}$.

\item We have the Hodge-Riemann bilinear relations:
\begin{enumerate}
\item $Q(B^{p, q}, B^{s, t})=0$ if $s\neq q$.

\item $(\sqrt{-1})^{-r}(-1)^qQ(\xi, \bar{\xi})>0$ if $0\neq \xi \in
B^{p, q}$ and $p+q=r$
\end{enumerate}
where $$Q(\tau,
\eta)=(-1)^{\frac{r(r+1)}{2}}\int_X\mathcal{L}^{n-r}(\tau\wedge
\eta)$$ and $\tau, \eta\in B^r$.
\end{enumerate}
\end{proposition}

\begin{proof}
By the assumption of the morphic conjecture, $\lambda$ is the
restriction of $\Lambda$ on $\mathscr{A}$. From the relation
$h=[\Lambda, L]$, we see that $h$ restricts to an operator on
$\mathscr{A}$. Thus $\mathcal{L}, \lambda, h$ give a
sub-$sl_2(\C)$-structure on $\mathscr{A}$ and therefore it admits a
sub-Lefschetz decomposition of the Lefschetz decomposition of
$H^*(X; \C)$ which is compatible with the sub-Hodge structure. The
Hard Lefschetz theorem is a formal consequence of the Lefschetz
decomposition (see e.g. \cite{Lew}, Chapter 11). The restriction of
the classical Hodge-Riemann bilinear relations to $\mathscr{A}$
gives the similar relations.
\end{proof}

We observe that the morphic signatures are independent of the odd
part $OH^b(X; \C)$ of the cohomology groups. In the following, by
assuming the morphic conjecture on $EH^a(X; \C)$, we are going to
generalize the classical Hodge index theorem.

\begin{theorem}(Hodge index theorem)
If the morphic conjecture is true on $EH^a(X; \C)$ where $X$ is a
connected projective manifold of dimension $m=2n$ and $a$ is a
nonnegative integer, then
$$\sigma_{a+n}(X)=\sum_{0\leq p, q\leq m}(-1)^qh^{p, q}_{\gamma(a, p,
q)}$$ where $\gamma(a, p, q)=\gamma(EH^a(X; \C), p, q)$.
\end{theorem}

\begin{proof}
Let $\omega$ be the $(1, 1)$-form associated to the standard
K\"{a}hler metric on $X$ and let $\mathscr{H}^{p, q}$ be the space
of harmonic $(p, q)$-forms on $X$. The Lefschetz operator $L:
\mathscr{H}^{p, q} \longrightarrow \mathscr{H}^{p+1, q+1}$ is
defined by $L\alpha=\omega \wedge \alpha$.

Since $\omega$ is integral, $L$ is an operator on $EH^a(X; \Q)$. We
have the Hodge inner product $< , >$ on $EH^a(X; \C)$, the adjoint
operator $\Lambda$ of $L$, the Hodge star operator $\bar{*}:
\mathscr{H}^{p, q} \longrightarrow \mathscr{H}^{n-p, n-q}$, and the
cup product pairing $(\ , \ )$ on $H^m(X; \C)$ satisfying $(\alpha,
\beta)=<\alpha, \overline{*}\beta>$ for all $\alpha, \beta\in H^m(X;
\C)$.

By the sub-Hodge structure on $\widetilde{L}^sH^k(X; \C)$, we
decompose $\widetilde{L}^sH^k(X;
\C)=\bigoplus_{p+q=k}\mathscr{H}^{p, q}_s$. Now we assume that the
morphic conjecture is true on $EH^a(X; \C)$. Then $EH^a(X; \C)$ is a
Hodge-Lefschetz algebra. The result now follows from Theorem
\ref{Abstract Hodge index theorem}.
\end{proof}

\begin{corollary}
When $a=n$, $h^{p, q}_{\gamma(n, p, q)}=h^{p, q}$, the above formula
gives the classical Hodge index theorem:
$$\sigma(X)=\sigma_{2n}(X)=\sum_{0\leq p, q\leq m}(-1)^qh^{p, q}_{\gamma(n, p, q)}=\sum_{0\leq p, q\leq
m}(-1)^qh^{p, q}.$$
\end{corollary}

\begin{example}
\begin{enumerate}
\item Suppose that $X$ is a complex projective surface. Then
$\sigma_1(X)=2-h^{1,1}(X),
\sigma_2(X)=\sigma(X)=2+2h^{2,0}(X)-h^{1,1}(X)$.

\item Suppose that $X$ is a general polarized abelian variety of dimension $g=2n$.
By a theorem of Mattuck, we have $H^{p, p}(X; \Q)\simeq \Q$ for
$0\leq p\leq g$ (see \cite{CH}, pg 559). Thus $\sigma_n(X)=1$. But
as a smooth manifold, $X$ is the boundary of a solid torus, hence
the signature of $X$ is 0.

\item
For a smooth hypersurface $X$ of dimension $m=2n$ in $\P^{2n+1}$,
for $p\neq n$, $H^{p, p}(X; \Q)$ is 1-dimensional and is generated
by algebraic cycles. Therefore the adjoint operator $$\Lambda :
H^{2n-p, 2n-p}(X; \Q) \longrightarrow H^{p, p}(X; \Q)$$ is an
isomorphism for $0\leq p<n$. For $p=n$, $\Lambda: H^{n, n}(X; \Q)
\longrightarrow H^{n, n}(X; \Q)$ is an isomorphism. The GSCA is
trivially true for this case, hence the signature formula
$$\sigma_{a+n}(X)=1+(-1)^{n-1}+(-1)^nh^{n, n}_{a+n}$$
is valid. In particular, we are especially interested in
$$\sigma_n(X)=1+(-1)^{n-1}+(-1)^nh^{n, n}_n$$
where $h^{n, n}_n$ is the dimension of the subspace of $H^{n, n}(X)$
which is generated by algebraic cycles. Thus any way to calculate
$\sigma_n(X)$ is equivalent to the calculation of $h^{n, n}_n$. The
Hodge conjecture predicts that $h^{n, n}_n=h^{n, n}_{\Q}$ where
$h^{n, n}_{\Q}$ is the dimension of $H^{n, n}(X; \Q):= H^{n, n}(X)
\cap H^{2n}(X; \Q)$.  We note that even for smooth hypersurfaces of
even dimension, the Hodge conjecture is known only for some small
degree.

Let $\sigma_n(X)=a-b$ the difference between the numbers of positive
and negative eigenvalues of $(\ , \ )$. By the Poincar\'e duality
theorem, $(\ , \ )$ is non-degenerate, hence $a+b=h^{n, n}_n$ and we
get $a=\frac{(-1)^n+1}{2}h^{n, n}_n$. Therefore, the cup product
pairing $(\ , \ )$ is positive definite on $H^{n, n}_n$ if $n$ is
even and negative definite on $H^{n, n}_n$ if $n$ is odd.
\end{enumerate}
\end{example}

\begin{acknowledgement}
The author thanks Blaine Lawson for his guidance for the author's
graduate study in Stony Brook, Christian Haesemeyer for his nice
comments and the secretaries in the Applied Mathematics Department
of National Chiao Tung University of Taiwan for their warmth. He
thanks Uwe Jannsen and the referee for their very nice comments and
corrections, the National Center of Theoretical Sciences of Taiwan
(Hsinchu) for providing a wonderful working environment.
\end{acknowledgement}

\bibliographystyle{amsplain}

\end{document}